\documentclass[11pt]{article}
\usepackage{color,amsmath,amssymb,amsfonts,amstext,amsthm,latexsym}

\usepackage{epsfig}
\usepackage{graphicx}
\usepackage{color}

  \newtheorem{Main Results}[theorem]{MainResults}

\begin{document}

\title{An alternative expression of Di Paola and Falson's formula for stochastic dynamics
  }

 \author{Xu Sun\\
 School of Mathematics and Statistics\\Huazhong University of Science and Technology\\ Wuhan 430074, Hubei, China \\
   E-mail: xsun@mail.hust.edu.cn; xsun15@gmail.com\\ \\
Jinqiao Duan\\
Institute for Pure and Applied Mathematics, University of California\\
Los Angeles, CA 90095\\ \& \\
Department of Applied Mathematics, Illinois Institute of Technology \\
   Chicago, IL 60616, USA \\
   E-mail: jduan@ipam.ucla.edu; duan@iit.edu\\
   \\
    Xiaofan Li \\
Department of Applied Mathematics, Illinois Institute of Technology \\   Chicago, IL 60616, USA \\   E-mail:    lix@iit.edu\\ \\
   }

\date{February 3, 2012}
\maketitle

\newpage
\begin{abstract}
Di Paola and Falsone's formula is widely used in studying stochastic dynamics of nonlinear systems under Poisson white noise. In this short communication, an alternative expression  is presented. Compared to Di Paola and Falsone's  original expression, the alternative one is applicable under more general condition, and shows significantly improved performance in numerical implementation. The alternative expression turns out to be a special case of the Marcus integrals.

\bigskip



\end{abstract}

\section {Introduction}

In the presence of random effects, dynamics of systems in engineering and science is often modeled by stochastic differential equatinos (SDEs) driven by Gaussian or non-Gaussian processes. The solutions of SDEs are often defined by  Ito or Stratonovich integrals. When the excitation is parametric Gaussian white noise, Stratonovich SDEs are different from the corresponding Ito SDEs: they have an additional term, the so-called WZ correction term \cite{CWSTo2000, WongZakai1965}. The studies in \cite{Ibrahim1985, YongLin1987} attribute the extra term to the conversion from physical white noise to ideal white noise. For Poisson white noise excitation, the correction terms for converting Stratonovich SDEs into Ito SDEs are first presented by Di Paola and Falsone \cite{DiPaolaFalsone1993, DiPaolaFalsone1993b}. Although some disagreement has been brought up concerning Di Paola and Falsone's correction term in earlier days \cite{Hu1994, Grigoriu1998}, the Stratonovich SDEs proposed in \cite{DiPaolaFalsone1993, DiPaolaFalsone1993b} have been verified extensively and are widely used to study stochastic dynamical systems under Poisson white noise excitation \cite{ZengZhu2010a, ZengZhu2010b, Proppe2002,CaddemiDiPaola1996}.

The main objective of this short communication is to give an alternative expression of Di Paola and Falsone's formula. Compared with the original one, the alternative formula is applicable under more general condition and shows better performance in numerical computation. Note that we  consider only the one-dimensional cases. The conclusion can be generalized to high-dimensional cases which, however,  will not be discussed in this short communication.

\section{Review of Di Paola and Falsone's formula}
Consider a dynamical system governed by
\begin{align}\label{section1_tmp1}
\dot Z(t)  =f(Z(t), t) + g(Z(t), t) \dot C(t),
\end{align}
where $Z$ represents the state variable, $f$ and $g$ are deterministic functions of $Z$ and $t$, $\dot C(t)$ is Poisson white noise expressed as
\begin{align}\label{section1_poissonwhitenoise}
\dot C(t) = \sum_{k=1}^{N(t)} R_k\delta(t-t_k).
\end{align}
In Eq. (\ref{section1_poissonwhitenoise}), $N(t)$ is a Poisson process with intensity parameter $\lambda$, $\delta(t-t_k)$ is a Dirac function at $t_k$, $R_k$ is a random variable representing the $k$-th impulse. The Poisson white noise $\dot C(t)$ is regarded as the formal derivative of some compound Poisson process $C(t)$,
\begin{align}\label{section1_tmpp1}
C(t)=\sum_{k=1}^{N(t)} R_k U(t-t_k),
\end{align}
where $U$ is a the step function.

In Ito sense, (\ref{section1_tmp1}) is often written in the form
\begin{align}\label{section1_sdei}
{\rm d} Z(t)  =f(Z(t), t)\,{\rm d}t  + g(Z(t), t)\, {\rm d} C(t),
\end{align}
Di Paola and Falsone suggest to interpret (\ref{section1_tmp1}) as \cite{DiPaolaFalsone1993, DiPaolaFalsone1993b}
\begin{align}\label{section1_sdes}
{\rm d} Z(t)  =f(Z(t-), t)\,{\rm d}t   + \sum_{j=1}^{\infty} \frac{g^{(j)}(Z(t-), t)}{j!}\, \left({\rm d}C(t)\right)^j,
\end{align}
where
\begin{align}\label{section1_dg}
\begin{cases}
g^{(1)} (x,t)=g(x,t),\\
g^{(j)}(x,t)= g(x,t) \frac{\partial}{\partial x} g^{(j-1)}(x,t), \quad \text{for $j\ge 2$.}
\end{cases}
\end{align}
(\ref{section1_sdes}) is often called the Stratonovich SDE for (\ref{section1_tmp1}) \cite{DiPaolaFalsone1993, DiPaolaFalsone1993b, ZengZhu2010a, ZengZhu2010b}.
 Note that ${\rm d}C(t)=C(t)-C(t-)$, and (\ref{section1_tmpp1}) implies that
 \begin{align}
 {\rm d}C(t)=\begin{cases}
 0, &\text{for}\quad t\ne t_k;\\
 R_k, &\text{for}\quad t=t_k.
 \end{cases}
 \end{align}

The effectiveness of  Di Paola and Falsone's formula expressed by (\ref{section1_sdes}) has been verified extensively \cite{ZengZhu2010a, ZengZhu2010b, Proppe2002,CaddemiDiPaola1996}.

\section{Derivation of the alternative formula}
Define $Y(\lambda, t)$ as
\begin{align}\label{section2_1}
Y(\lambda, t)=\sum_{j=1}^{\infty} \frac{g^{(j)}(Z(t-), t)}{j!}  \lambda  ^j,
\end{align}
where $g^{(j)}(x,t)$ ($j=1, 2, \cdots$) are defined by (\ref{section1_dg}). Using the notation $Y$ defined in (\ref{section2_1}), we can write the Stratonovich SDE (\ref{section1_sdes}) as
\begin{align}\label{section2_2}
{\rm d} Z(t)  =f(Z(t), t)\,{\rm d}t  + Y({\rm d}C(t),t).
\end{align}
Construct an ordinary differential equation(ODE) as
\begin{align}\label{section2_ode}
\begin{cases}
\frac{{\rm d} y(\lambda)}{{\rm d} \lambda}=g(Z(t-)+y(\lambda), t) \quad \text{for} \quad  0\le\lambda\le {\rm d}\,C(t), \\
y(0)=0.
\end{cases}
\end{align}
We claim that $Y({\rm d}\,C(t),t)$ as defined in \eqref{section2_1} is the solution of the above ODE~(\ref{section2_ode}) at $\lambda={\rm d}C(t)$. The proof is given as follows.

First, we check  the existence and uniqueness of solution for (\ref{section2_ode}). It is well known that if the function $g(x,t)$ is Lipschitz continuous, then the existence and uniqueness of (\ref{section2_ode}) can be guaranteed. Since (\ref{section2_1}) requires the function $g(x, t)$ be analytic with respect to $x$ \cite{DiPaolaFalsone1993, DiPaolaFalsone1993b}, which implies Lipschitz continuity, the existence and uniqueness of the solution follows.

Take Taylor expansion of $g(Z(t-)+y(\lambda, t), t)$ at $\lambda=0$, then (\ref{section2_ode}) becomes
\begin{align}\label{section2_4}
\frac{{\rm d}\, y(\lambda)}{{\rm d}\, \lambda}= g(z(\lambda ), t)\big|_{\lambda=0}  + \sum_{k=1}^{\infty} \frac{\partial ^k}{\partial \lambda^k} g(z(\lambda ), t)\bigg|_{\lambda=0} \frac{ \lambda^{k}}{k!},
\end{align}
where
\begin{align}\label{section2_3}
z(\lambda )=Z(t-)+y(\lambda).
\end{align}

In the following, we  show by induction that
 \begin{align}\label{section2_5}
  \frac{\partial ^k}{\partial \lambda^k} g(z(\lambda ), t)   =g^{(k+1)} (z(\lambda ), t), \quad \text{for}\quad k\ge 1.
 \end{align}
For $k=1$, we obtain from (\ref{section2_3}) and (\ref{section2_ode}) that
 \begin{align} \label{section2_31}
\frac{\partial }{\partial \lambda } g(z(\lambda),t) & =\left(\frac{{\rm d} }{{\rm d}\,\lambda } z(\lambda )\right) \left( \frac{\partial }{\partial z (\lambda)} g(z(\lambda), t)) \right)\nonumber \\
& =\left(\frac{{\rm d} }{{\rm d}\,\lambda } y(\lambda )\right) \left( \frac{\partial }{\partial z (\lambda)} g(z(\lambda), t)) \right) \nonumber \\
& = g(z(\lambda), t)) \left( \frac{\partial }{\partial z (\lambda)} g(z(\lambda), t)) \right) \nonumber \\
& = g^{(2)} (z(\lambda), t)),
\end{align}
which means that (\ref{section2_5}) is true for $k=1$.

Suppose  (\ref{section2_5}) is true for $k=n$, i.e.,
  \begin{align} \label{section2_3n}
\frac{\partial^n }{\partial \lambda^n } g(z(\lambda),t) = g^{(n+1)} (z(\lambda), t)),
\end{align}
 then we get
\begin{align}
\frac{\partial^{n+1 }}{\partial \lambda^{n+1 }} g(z(\lambda),t)&= \frac{\partial}{\partial \lambda}   \left(\frac{\partial^{n}}{\partial  \lambda^n } g(z(\lambda, t))\right) \nonumber\\
&=\frac{\partial}{\partial \lambda}    g^{(n+1)} (z(\lambda), t)) \nonumber\\
& =\left(\frac{{\rm d} }{{\rm d}\,\lambda } z(\lambda )\right) \left( \frac{\partial }{\partial z (\lambda)} g^{(n+1)}(z(\lambda), t)) \right)\nonumber \\
& = g(z(\lambda), t)) \left( \frac{\partial }{\partial z (\lambda)} g^{(n+1)}(z(\lambda), t)) \right) \nonumber \\
& = g^{(n+2)} (z(\lambda), t)),
\end{align}
which means that (\ref{section2_5}) is also true for $k=n+1$. Now the proof of (\ref{section2_5}) is complete.

Substituting (\ref{section2_5}) into (\ref{section2_4}), and using the fact that
\begin{align}
g^{(n+1)} (z(\lambda), t))\big|_{\lambda=0} = g^{(n+1)} (Z(t-), t)),
\end{align}
we get
\begin{align} \label{section2_6}
\frac{{\rm d}\, y(\lambda)}{{\rm d}\, \lambda}=   \sum_{k=1}^{\infty} g^{(k)}(Z(t-), t) \frac{ \lambda^{k-1}}{(k-1)!}.
\end{align}
Integrating from $0$ to $\lambda$ at both sides of (\ref{section2_6}), we get
\begin{align}\label{section2_7}
 y(\lambda) =   \sum_{k=1}^{\infty} g^{(k)}(Z(t-), t) \frac{ \lambda^{k }}{k!}\;.
\end{align}
Comparing (\ref{section2_7}) with (\ref{section2_1}), we see the equivalence
\begin{align} \label{equiv}
  Y({\rm d}C(t), t) =  y({\rm d}C(t)),
\end{align}
where the right-hand side $y({\rm d}C(t))$ is the solution of the initial value problem (\ref{section2_ode}) at $\lambda ={\rm d}C(t)$. In other words,
we obtain an equivalent and alternative expression for the Stratonovich SDE
(\ref{section1_sdes}) or \eqref{section2_2}
\begin{align}\label{ssde}
{\rm d} Z(t)  =f(Z(t), t)\,{\rm d}t  + y({\rm d}C(t)),
\end{align}
where $y({\rm d}C(t))$ is the solution of the ODE (\ref{section2_ode}) at $\lambda ={\rm d}C(t)$.

\section{Some remark and examples}
As shown in section $3$,  when $g(x,t)$ is analytic, $Y({\rm d}C(t), t)$
defined in (\ref{section2_1}) is equivalent to $y({\rm d}C(t))$ obtained from the ODE (\ref{section2_ode}), i.e., the two forms of SDEs, \eqref{section1_sdes} and \eqref{ssde} are equivalent.  However, compared with (\ref{section1_sdes}), (\ref{ssde}) has significant advantages in the following aspects.

First, for a given $g(x, t)$, one immediate advantage of the alternative expression \eqref{ssde} over (\ref{section1_sdes}) is that it is much easier to obtain analytical solution from the ODE \eqref{section2_ode} than finding the sum of the infinite Taylor series in \eqref{section1_sdes}. Second, the proposed alternative expression (\ref{ssde}) is applicable under much more general condition than the original one (\ref{section1_sdes}). For example, when $g(x,t)$ is only continuous or piecewise continuous rather than analytic, (\ref{ssde}) is still valid while (\ref{section1_sdes}) fails to exist. More importantly, the new expression (\ref{ssde}) could have significantly improved performance in numerical implementation than the original expression (\ref{section1_sdes}). The form in (\ref{section2_1}) is actually the Taylor expansion of the solution of (\ref{section2_ode}), which requires computing the high-order derivatives of the function $g(x,t)$ and possibly take a large number of terms to achieve high accuracy. However, since the alternative form (\ref{ssde}) is given in form of ODE, it can be solved by numerous algorithms such as Runge-Kutta or multistep methods, which could have much better convergence and stability properties than Taylor expansion methods.  The readers are referred to any numerical analysis book such as \cite{KincaidCheney2002} for more details about the Runge-Kutta, multistep, or Taylor expansion methods for solving ordinary differential equations numerically.

We note that the alternative expression, as given in  (\ref{ssde}), turns out to be a special case of the Marcus integrals, which have been defined for general semimartingale with jumps. For more details about the Marcus integral and its application, readers are referred to \cite{Marcus1978, Marcus1981, KurtzPardouxProtter1995, Applebaum2009}.

 {{\bf Example }}\\
  In this example, $f(x,t)$ and $g(x,t)x$ are chosen such that the analytical solution of the SDE is available to compare with the numerical results.

 Take $f(x,t)=x$, and $g(x,t)=x$, the SDE (\ref{section1_sdes}) becomes
 \begin{align}\label{example_5}
 {\rm d}Z(t) = Z(t){\rm d} t +Y({\rm d} C(t)),
 \end{align}
where $Y({\rm d} C(t))$ is expressed, by Di Paola and Falsone's formula, as
\begin{align}
 Y({\rm d} C(t), t)) =Z(t-) \sum_{k=1}^{\infty} \frac{\left({\rm d} C(t)\right)^k}{k!},
 \end{align}
or by the alternative expression, as the solution of the following ODE at $\lambda ={\rm d} C(t)$,
 \begin{align}\label{example_6}
\begin{cases}
\frac{{\rm d} y(\lambda)}{{\rm d} \lambda}= Z(t-)+y(\lambda)   \quad \text{for} \quad  0\le\lambda\le {\rm d}\,C(t), \\
y(0)=0.
\end{cases}
\end{align}

It is easy to check that the theoretical solution of (\ref{example_5}) is $e^{t+C(t)}$.

In the simulation, the compound Poisson process $C(t)$ is taken as
\begin{align}
C(t)=\sum_{k=1}^{N(t)} R_k U(t-t_k).
\end{align}
wher $N(t)$ is a Poisson process with intensity parameter as $10$ and $ R_k$ ($k=1, \cdots, N(t)$) are normally distributed random numbers independent of each other.

When no jumps occur, $  Y({\rm d} C(t))=0$ and
We solve  (\ref{example_5}) by second order Runge-Kutta method with time step as  $0.01$. We compute the jump of the solution $Z(t)$ by Di Paola and Falsone's formula and the alternative formula, respectively. For the former, we truncate the infinite consequence by taking the first $6$ terms, and for the latter,  we solve ODE (\ref{example_6}) numerically by uisng second order Runge-Kuta method with step size equal to $0.1$.  Figure 1 shows one random path of the Poisson compound process, and Figure 2 compares the corresponding computational results with the theoretical one. It can be seen from Fig. 2 that when the jumps are computed by the alternative formula, the numerical result  agrees very well with the theoretical one. However, when using Di Paola and Falsone's formula, the numerical result deviates from the theoretical one significantly, especially when the jump size is relatively large. We note that the accuracy for Di Paola and Falsone's formula could be improved by keeping more terms when truncating the infinite series. However, the number of terms needed to achieve the desired accuracy with Di Paola and Falsone's formula is not easy to predict. The alternative formula is preferred partly because its accuracy can be easily controlled by the step size and the order the Runge-Kutta or other methods.
\newpage
\begin{figure}[htp]
\vspace{4cm}
  \epsfig{file=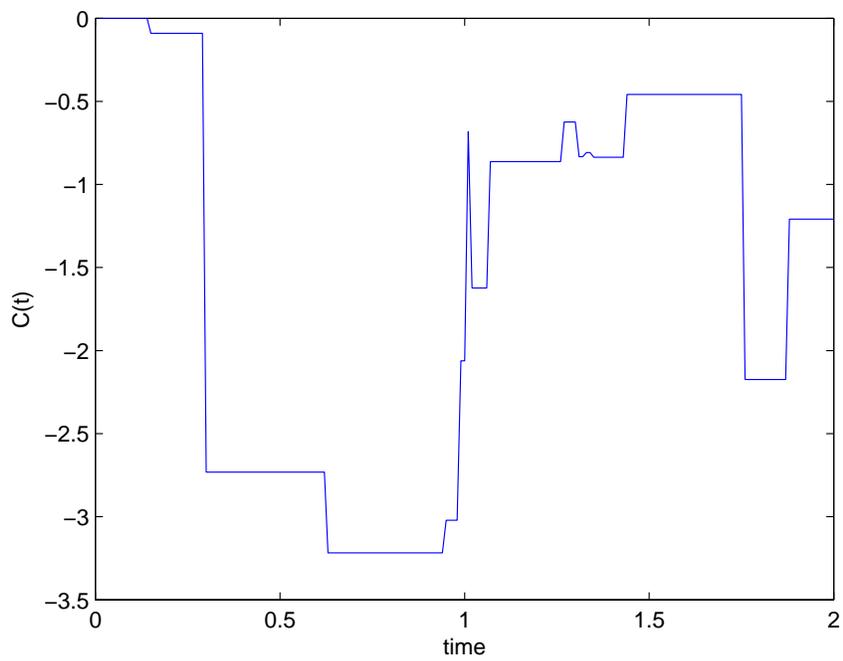,width=\linewidth}
  \caption{One random path of the driving process }
\end{figure}
\newpage
\begin{figure}[htp]
\vspace{4cm}
  \epsfig{file=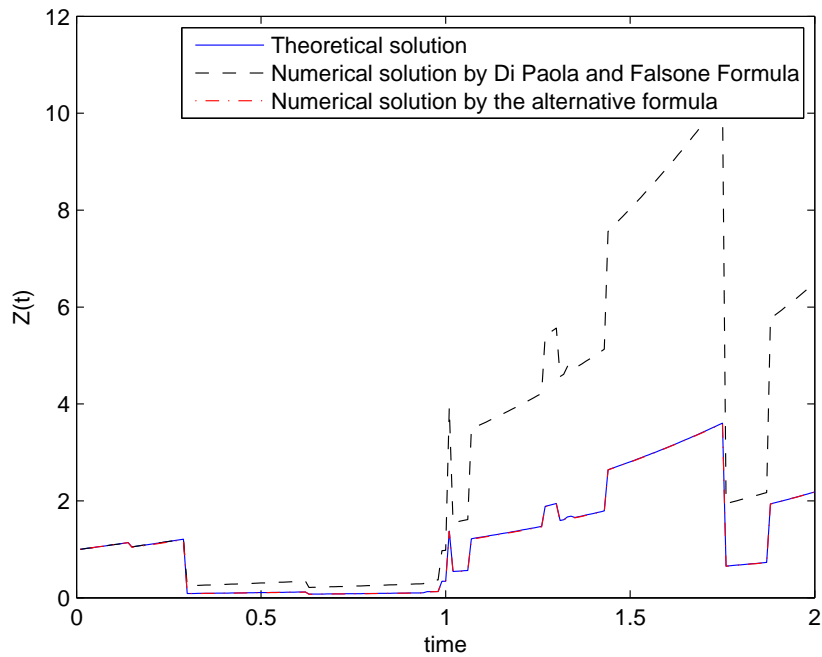,width=\linewidth}
  \caption{Comparison of the numerical results with the theoretical one}
\end{figure}
%

\end{document}